\documentclass[12pt]{amsart}

\usepackage{amsmath,graphics,comment}
\usepackage{amsfonts,amssymb}

  \includecomment{inlong} \excludecomment{inshort}   
 \includecomment{inlow} \excludecomment{inhigh}

\theoremstyle{plain}
\newtheorem*{theorem*}{Theorem}
\newtheorem*{lemma*} {Lemma}
\newtheorem*{corollary*} {Corollary}
\newtheorem*{proposition*} {Proposition}
\newtheorem{theorem}{Theorem}[section]
\newtheorem{lemma}[theorem]{Lemma}
\newtheorem{corollary}[theorem]{Corollary}
\newtheorem{proposition}[theorem]{Proposition}

\theoremstyle{remark}
\newtheorem*{remark}{Remark}
\newtheorem*{definition}{Definition}

\newtheorem*{claim}{Claim}

\theoremstyle{definition}

\def \R {\mathbf{R}}
\def \Z {\mathbf{Z}}
\def \C {\mathbf{C}}

\def \L {\mathbf{L}}

\def\eps{\epsilon}

\def\s{\sigma}

\def\Q{\Bbb{Q}}

\def\sign{\mbox{sign}}
\def\Z{\Bbb{Z}}
\def\C{\Bbb{C}}

\def\N{\Bbb{N}}
\def\l{\lambda}
\def\part{\partial}

\def\a{\alpha}
\def\g{\gamma}

\def\bp{\begin{pmatrix}}
\def\Arf{\mbox{Arf}}
\def\sm{\setminus}
\def\ep{\end{pmatrix}}
\def\bn{\begin{enumerate}}
\def\Rightleftarrow{\Leftrightarrow}

\def\rk{\mbox{rank}}
\def\en{\end{enumerate}}
\def\ba{\begin{array}}
\def\ea{\end{array}}

\def\L{\Lambda}
\def\S{\Sigma}
\def\s{\sigma}
\def\a{\alpha}
\def\b{\beta}
\def\ti{\tilde}

\def\fr12{\frac{1}{2}}

\def\COT{Cochran, Orr and Teichner }

\def\im{\mbox{Im}}

\def\ker{\mbox{Ker}}

\def\v{\varphi}

\def\lteta{\eta^{(2)}}

\def\varnothing{\emptyset}

\def\cmtbf#1{} \def\cmt#1{}

\begin{document}

\title{$L^2$--eta--invariants and their approximation by unitary eta--invariants}
\author{Stefan Friedl}
\date{\today}
\begin{abstract}
Cochran, Orr and Teichner introduced $L^2$--eta--invariants to detect highly non--trivial examples of non slice knots.
Using a recent theorem by L\"uck and Schick we show that their metabelian $L^2$--eta--invariants can be viewed
as the limit of finite dimensional unitary representations. We recall a ribbon obstruction theorem proved by the author
using finite dimensional unitary eta--invariants. We show that if for a knot $K$ this ribbon obstruction vanishes then
the  metabelian $L^2$--eta--invariant vanishes too. The converse has been shown by the author not to be true.
\end{abstract}
\maketitle


\section{Introduction}
 
A knot $K\subset S^{n+2}$ is a smooth submanifold homeomorphic to $S^n$. A knot is called slice if it bounds a smooth disk in
$D^4$.  We say that a knot $K$ is algebraically  slice if $K$ has a Seifert matrix of the form
$\bp 0 & B \\ C&D \ep$ where $B,C,D$ are square matrices of the same size. It is a well-known fact that
any slice knot is algebraically slice.
Levine showed that  in higher odd
dimensions the converse is true, i.e. if a knot is algebraically slice it is also geometrically slice (cf. \cite{L69}).
In the classical dimension $n=1$ this no longer holds as was shown by Cassan and Gordon \cite{CG86}.

A knot
$K\subset S^3$ is called ribbon if there exists a smooth 
disk $D$ in $S^3\times [0,1] \subset D^4$ ($S^3=S\times 0$)
bounding $K$ such that the projection map $S^3\times [0,1]\to [0,1]$ is a Morse map and has no local minima.
Such a slice disk is called a ribbon disk.
Fox \cite{F61} conjectured that all slice knots are ribbon. 
 
In   \cite{F03} the author studies metabelian unitary eta--invariants of $M_K$, the result of 
zero framed surgery along a knot $K\subset S^3$. These can be used  
  to detect knots which are not slice respectively, not
ribbon.

For a pair $(M^3,\varphi:\pi_1(M)\to G)$ Cheeger and Gromov \cite{CG85} introduced the $L^2$--eta--invariant
$\eta^{(2)}(M,\varphi )$.
 Cochran, Orr and Teichner \cite{COT01} gave examples of knots which look slice 
`up to a certain level' but can be shown to be not slice using $L^2$--eta--invariants.  

L\"uck and Schick \cite{LS01} showed that $L^2$--eta--invariants can be viewed as a limit of ordinary unitary
eta--invariants if $G$ is residually finite.
We show that the metabelian  groups used by Cochran, Orr and Teichner
are  residually finite. Sorting out several technical problems we can show that if for a knot $K$
the metabelian eta--invariant ribbonness obstruction vanishes then the metabelian $L^2$--eta--invariant sliceness obstruction
vanishes as well. In \cite{F03} we show that the converse is not true.

The structure of the paper is as follows. In section \ref{sectioneta} we recall the   eta--invariant
sliceness and ribbonness obstruction theorems of \cite{F03}.  In section \ref{sectioncot} we give the definition of
$(n)$--solvability for a knot
$n\in \frac{1}{2}\N$, and quote some results of \cite{COT01}.
Furthermore we state the metabelian $L^2$--eta--invariant sliceness obstruction theorem of Cochran, Orr and Teichner.
 We state and prove the main theorem in section \ref{sectionmainthm}. 

{ \bf Acknowledgment.} I would like to thank Jerry Levine, Kent Orr and Taehee Kim for helpful discussions
and comments.
 
\section{Unitary eta--invariants as knot invariants} \label{sectioneta}
Let $M^{2q+1}$ be a closed odd-dimensional smooth manifold and $\a:\pi_1(M) \to U(k)$ a unitary representation.
Atiyah, Patodi, Singer \cite{APS75} associated to $(M,\a)$ a number $\eta(M,\a)$ called the (reduced)
eta--invariant of
$(M,\a)$. This invariant has the property that if  
 $\partial(W^{2q+2},\b)=(M^{2q+1},\a)$ then
\[ \eta(M,\a) =\sign_{\b}(W)-k\sign(W) \]
where $\sign_{\b}(W)$ denotes the signature of $W$ twisted by $\b$.

\subsection{Abelian eta--invariants}
 Let $K$ be knot, $\mu$ a meridian and $A$ a Seifert matrix for $K$.
Let $\a:\pi_1(M_K)\to U(1)$ be a representation, then 
\[ \eta(M_K,\a)=\s_{z}(K):=\sign(A(1-z)+A^t(1-\bar{z}))\]
where $z:=\a(\mu)$  (cf. \cite{L84}).

The following proposition follows immediately from the definitions and the explicit computation of the abelian
eta--invariant. 

\begin{proposition}
Let $K$ be an algebraically slice knot, then $\eta(M_K,\a)=0$ for any representation
$\a:\pi_1(M_K)\to U(1)$ which sends the meridian to a transcendental number.
\end{proposition}

If a knot satisfies the conclusion of this proposition we say that $K$ has zero abelian eta--invariant sliceness obstruction.
\subsection{Metabelian eta--invariants}
There exists a canonical map $\eps: \pi_1(M_K)\to H_1(M_K)=\Z$ sending the meridian to $1$.
Denote the $k$--fold cover of $M_K$ by $M_k$.
If $k$ is a prime power, then Casson and Gordon \cite{CG86} showed that
$H_1(M_k)=\Z \oplus TH_1(M_k)$ where $TH_1(M_k)$ denotes
the  $\Z$--torsion part of $H_1(M_k)$.
Furthermore there exists a non--singular symmetric linking pairing
\[ \l_{lk}:TH_1(M_k)\times TH_1(M_k)\to \Q/\Z \]
We say that $P_k \subset TH_1(M_k)$ is a $\L$--metabolizer for $\l_{lk}$ if $P_k$ is a $\L$--submodule and if
\[ P_k=P_k^{\perp}:=\{ x\in TH_1(M_k) | \l_{lk}(x,y)=0 \mbox{ for all }y\in TH_1(M_k)\} \]

Denote by $\ti{M}_K$ the universal abelian cover corresponding to $\eps$. $H_1(\ti{M}_K)$ carries
a $\L:=\Z[t,t^{-1}]$--module structure, we will henceforth write  $H_1(M_K,\L)$ for $H_1(\ti{M}_K)$.  
Blanchfield \cite{B57} shows that there exists a non-singular $\L$-hermitian pairing
\[ \l_{Bl}: H_1(M_K,\L)\times H_1(M_K,\L) \to \Q(t)/\L\]
For a $\L$-submodule $P\subset H_1(M_K,\L)$ define  
\[ P^{\perp}:=\{ v\in H_1(M_K,\L) | \l_{Bl}(v,w)=0 \mbox{ for all } w \in P\}
\]  
If $P \subset H_1(M_K,\L)$ is such that $P=P^{\perp}$, then we say  
that $P$ is a metabolizer for $\l_{Bl}$ and that $\l_{Bl}$ is metabolic. 
Note that Kearton \cite{K75} showed that a knot is algebraically slice if and only if
$\l_{Bl}$ is metabolic.
\\

Recall that for a group $G$ the central series is defined inductively by $G^{(0)}:=G$ and $G^{(i)}:=[G^{(i-1)},G^{(i-1)}]$.
Let $\pi:=\pi_1(M_K)$. We study metabelian representations,
i.e. representations that factor through $\pi/\pi^{(2)}$. Consider
\[ 1\to \pi^{(1)}/\pi^{(2)}\to \pi/\pi^{(2)}\to\pi/\pi^{(1)}\to 1 \]
Note that $\pi^{(1)}/\pi^{(2)}\cong H_1(\ti{M}_K)$ and $\pi/\pi^{(1)}= H_1(M_K)=\Z$, in particular this sequence splits and we
get an isomorphism
\[   \pi/\pi^{(2)}
 \cong  \Z \ltimes H_1(M_K,\L)  \]
where $1 \in \Z$ acts by conjugating with $\mu$ respectively by multiplying
by $t$. Eta invariants corresponding to metabelian representations in the context of knot theory were first
studied by Letsche \cite{L00}.

For a group $G$ denote by
$R_k^{irr}(G)$ (resp. ${R}_k^{irr,met}(G)$) the
set of  irreducible, $k$-dimensional,  unitary
(metabelian)
representations of
$G$. By $\hat{R}$ we denote the conjugacy classes of such representations.
The above discussion shows that for a knot $K$ we can identify
\[ R_k^{irr,met}(\pi_1(M_K)) = R_k^{irr}(\Z \ltimes H_1(M_K,\L)) \]

\begin{lemma} \cite{F03} \label{lemma1}
Let $z \in S^1$ and $\chi:H_1(M,\L) \to H_1(M,\L)/(t^k-1) \to S^1$ a character. Then
 \[ \ba{rcl} \a_{(k,z,\chi)} \hspace{-0.05cm}= \hspace{-0.05cm}\a_{(z,\chi)} :
\Z  \hspace{-0.05cm}\ltimes  \hspace{-0.05cm}H_1(M,\L)& \hspace{-0.1cm}\to 
\hspace{-0.1cm} &U(k) \\
    (n,h) & \hspace{-0.1cm}\mapsto \hspace{-0.1cm} &
z^n  \hspace{-0.1cm} \bp 0& \dots &0&1 \\ 1&\dots &0&0 \\
\vdots &\ddots &&\vdots \\
     0&\dots &1&0 \ep^n  \hspace{-0.1cm} \bp \chi(h) &0&\dots &0 \\
 0&\chi(th) &\dots &0 \\
\vdots &&\ddots &\vdots \\ 0&0&\dots &\chi(t^{k-1}h) \ep \ea \]
defines a  representation.

Conversely any irreducible representation $\a \in{R}_k^{irr}(\Z \ltimes H_1(M,\L))$ is (unitary) conjugate to 
$\a_{(z,\chi)}$ for some $z\in S^1$ and a character $\chi:H_1(M,\L) \to H_1(M,\L)/(t^k-1) \to S^1$
which does not factor through $H_1(M,\L)/(t^l-1)$ for some $l<k$.
\end{lemma}

We denote by $P^{met}_k(\pi_1(M_K))$ the set of metabelian representations of $\pi_1(M_K)$
that are conjugate to $\a_{(z,\chi)}$  with $z$ transcendental and $\chi$ of prime power order.
Furthermore for $p$ a prime we write $ P_{k,p}^{irr,met}(\pi_1(M_K))$ for the set of representations
where $\chi$ has order a power of $p$.
In \cite{F03} we  prove the following sliceness obstruction theorem which is the strongest 
theorem detecting non--torsion knots which is not based on $L^2$--eta--invariants.

\begin{theorem} \label{mainthm2}
Let $K$ be a slice knot, $k_1,\dots,k_r$ pairwise coprime prime powers, then there
exist  $\L$--metabolizers $P_{k_i} \subset TH_1(M_{k_i}), i=1,\dots,r$ for the linking pairings $\l_{k_i}$, such that
for any prime number $p$ and any choice of  
irreducible representations 
$\a_i:\pi_1(M_K) \to \Z \ltimes H_1(M_K,\L)/(t^{k_i}-1) \to U(k)$ vanishing on $0
\times P_{k_i}$ and lying in
$ P_{k_i,p}^{irr,met}(\pi_1(M_K))$ we get $\eta(M_K,\a_1 \otimes \dots \otimes \a_r)=0$.
\end{theorem}

If a knot $K$ satisfies the conclusion of this theorem we say that $K$ has zero metabelian eta--invariant sliceness
obstruction.

In \cite{F03} we prove the following ribbon obstruction theorem.
In the proof we only use the well-known fact that if $K$ is ribbon then $K$ has a slice  disk $D$ such that 
$\pi_1(S^3 \sm K) \to  \pi_1(D^4\sm D)$ is surjective.

\begin{theorem}\cite{F03} \label{mainribbonthm}
Let $K \subset S^3$ be a ribbon knot.
Then there exists a metabolizer $P$ for the Blanchfield pairing such that
 for any $\a_{(z,\chi)}$ with $z$ transcendental and $\chi$ of prime power order, vanishing
on $0\times P$ we get $\eta(M_K,\a_{(z,\chi)})=0$.
\end{theorem}

We say that $K$ has zero metabelian eta--invariant ribbonness obstruction if the conclusion of the
theorem holds for $K$.

\cmt{----------------------------------------------------------------------------------}
\section{The Cochran--Orr--Teichner sliceness obstruction} \label{sectioncot}
\cmt{--------------------------------}
\subsection{The Cochran--Orr--Teichner  sliceness filtration}

We give a short introduction to the sliceness filtration introduced by Cochran, Orr and Teichner \cite{COT01}.
For a manifold $W$ denote by $W^{(n)}$ the cover corresponding to $\pi_1(W)^{(n)}$. Denote the
equivariant intersection form
\[ H_2(W^{(n)})\times H_2(W^{(n)}) \to \Z [\pi_1(W)/  \pi_1(W)^{(n)}] \]  
by $\lambda_n$, and the
self-intersection form by $\mu_n$. 
An $(n)$-Lagrangian is a submodule $L \subset H_2(W^{(n)})$ on which $\lambda_n$ and $\mu_n$ vanish
and which maps onto a Lagrangian of $\lambda_0 :H_2(W)\times H_2(W)\to \Z$.

\begin{definition}\cite[def. 8.5]{COT01}
A knot $K$ is called $(n)$-solvable if $M_K$ bounds a spin 4-manifold $W$ such that $H_1(M_K)\to H_1(W)$ is an
isomorphism and such that $W$ admits two dual $(n)$-Lagrangians. This means that $\l_n$ pairs the
two Lagrangians non-singularly and that the projections freely generate $H_2(W)$. 

A knot $K$ is called $(n.5)$-solvable if $M_K$ bounds a spin 4-manifold $W$ such that $H_1(M_K)\to H_1(W)$ is an
isomorphism and such that $W$ admits an $(n)$-Lagrangian and a dual $(n+1)$-Lagrangian.
 
We call $W$ an $(n)$-solution respectively $(n.5)$-solution for $K$. 
 \end{definition}

\begin{remark}
\bn
\item The size of an $(n)$-Lagrangian depends only on the size of $H_2(W)$,
in particular if $K$ is slice, $D$ a slice disk, then $\overline{D^4\sm N(D)}$
is an
$(n)$-solution for
$K$ for all
$n$, since
 $H_2(D^4\sm N(D))=0$.
\item By the naturality of covering spaces and homology with twisted coefficients it follows that
if $K$ is $(h)$-solvable, then it is $(k)$-solvable for all $k<h$.
\en
\end{remark}

\begin{theorem}
\[ \ba{rcl} K \mbox{ is (0)-solvable} & \Leftrightarrow & \Arf(K)=0 \\
K \mbox{ is (0.5)-solvable} & \Rightleftarrow & K \mbox{ is algebraically slice} \\
K \mbox{ is (1.5)-solvable} & \Rightarrow & \mbox{Casson-Gordon   invariants
vanish and $K$ algebraically slice} \ea \]
The converse of the last statement is not true, i.e. there exist algebraically slice knots which have
zero Casson-Gordon invariants but are not $(1.5)$-solvable.
\end{theorem}

The first part, the third part and the $\Leftarrow$ direction of the second part have been shown
by \COT  \cite[p. 6, p. 72, p. 66, p. 73]{COT01}.
\COT \cite[p. 6]{COT01} showed that a knot is $(0.5)$ solvable if and only if the
Cappell-Shaneson surgery obstruction in $\Gamma_0(\Z[\Z]\to \Z)$ vanishes. This is equivalent to a knot being algebraically
slice (cf. \cite{K89}).
  Taehee Kim \cite{K02} showed that there exist $(1.0)$-solvable knots which have zero Casson-Gordon invariants, but
 are not $(1.5)$-solvable. Cochran, Orr and Teichner \cite{COT01} also showed that there exist $(2)$--solvable knots
which are not $(2.5)$--solvable.

\subsection{$L^2$--eta--invariants as sliceness-obstructions}
In this section we'll very  quickly summarize some $L^2$--eta--invariant theory.

Let $M^3$ be a smooth manifold and $\varphi:\pi_1(M)\to G$ a homomorphism,
then Cheeger and Gromov \cite{CG85} defined an invariant
$\eta^{(2)}(M,\varphi)\in \R$, the (reduced)
$L^2$--eta--invariant. When it's clear which homomorphism we mean, we'll write 
$\eta^{(2)}(M,G)$ for $\eta^{(2)}(M,\varphi)$.

\begin{remark}
If  
$\partial(W,\psi)=(M^3,\varphi)$, then (cf. \cite[lemma 5.9 and remark 5.10]{COT01})
\[ \eta^{(2)}(M,\varphi)=\sign^{(2)}(W,\psi)-\sign(W) \]
where $\sign^{(2)}(W,\psi)$  denotes Atiyah's
   $L^2$-signature (cf. \cite{A76}).
\end{remark}

\COT study when $L^2$--eta--invariants vanish for homomorphisms $\pi_1(M_K)\to G$, where $G$ is a PTFA-group.
PTFA stands for poly-torsion-free-abelian, and means that there exists a normal subsequence where each quotient is
torsion-free-abelian.

\begin{theorem}{\cite[p. 5]{COT02}} \label{propl2etazero}
Let $G$ be a PTFA-group with $G^{(n)}=1$.
If $K$ is a knot, and $\varphi:\pi_1(M_K)\to G$ a homomorphism which extends over a $(n.5)$-solution of $M_K$, then
$\lteta(M_K,\varphi)=0$. In particular if
$K$ is slice and $\varphi$ extends over $D^4\sm D$ for some slice disk $D$, then
$\lteta(M_K,\varphi)=0$.
\end{theorem}

\begin{remark}
It's a crucial ingredient in the proposition 
that the group $G$ is a PTFA-group, for example it's not true in
general that $\lteta(M_K,\Z/k)=0$ for a slice knot $K$. Corollary \ref{corl2equalssum} shows that
$\lteta(M_K,\Z/k)=\sum_{j=1}^k \s_{e^{2\pi ij/k}}(K)$, 
but this can be non-zero 
for some slice knot $K$,  e.g. take a slice knot with Seifert matrix
\[ A =\bp 0&0&1&1 \\ 0&0&0&1 \\ 1&1&0&1 \\ 0&1&0&0\ep \]
Then $\lteta(M_K,\Z/6)=-2$.
\end{remark}

We use this theorem only in the abelian and the metabelian setting.
Let $\Q\L:=\Q[t,t^{-1}]$.

\begin{theorem}{\cite{COT01}} 
\bn
\item If $K$ is $(0.5)$--solvable, then $\eta^{(2)}(M_K,\Z)=0$.
\item
If $K$ is $(1.5)$--solvable, then
 there exists a metabolizer $P_{\Q}\subset H_1(M_K,\Q\L)$ for the rational Blanchfield pairing
\[ \l_{Bl,\Q}: H_1(M_K,\Q\L) \times H_1(M_K,\Q\L) \to \Q(t)/\Q[t,t^{-1}] \]
such that for all $x\in P_{\Q}$ we get $\lteta(M_K,\b_x)=0$ where $\b_x$ denotes the map
\[ \pi_1(M_K)\to \Z \ltimes H_1(M_K,\L) \to \Z \ltimes H_1(M_K,\Q\L)
\xrightarrow{id \times \l_{Bl,\Q}(x,-)}\Z \ltimes \Q(t)/\Q[t,t^{-1}]\]
\en
 \end{theorem} 

\begin{proof}
Let $D$ be a slice disk for $K$, write $N_D:=\overline{D^4\sm N(D)}$.
Then the statement follows from  proposition  \ref{propl2etazero} and work by Letsche \cite{L00} who
showed that for  
$P_{\Q}:=\ker\{H_1(M_K,\Q\L)\to H_1(N_D,\L\Q)\}$ the map $\b_x$ extends over
$\pi_1(N_D)$.  
\end{proof}

We say that $K$ has zero abelian  $L^2$--eta--invariant sliceness obstruction  
if  $\eta^{(2)}(M_K,\Z)=0$.
We say that $K$ has zero metabelian  $L^2$--eta--invariant sliceness obstruction  
if there exists a metabolizer $P_{\Q}\subset H_1(M_K,\Q\L)$ for $\l_{Bl,\Q}$
such that for all $x\in P_{\Q}$ we get $\lteta(M_K,\b_x)=0$.

%

\section{Relation between eta--invariants and $L^2$--eta--invariants} \label{sectionmainthm}
If a knot $K$ has zero abelian eta--invariant sliceness obstruction, then a multiple of $K$ is algebraically slice
(cf. Levine \cite{L69b} and Matumuto \cite{M77}), in particular $K$ has zero abelian $L^2$--eta--invariant
sliceness obstruction. This fact will also follow immediately from corollary \ref{corl2equalsint}.
Conversely, if $K$ has zero abelian $L^2$--eta--invariant,  then it is not necessarily true that 
$K$ has zero abelian eta--invariant, as was shown in \cite{F03}. 

In \cite{K02} Taehee Kim gave examples of knots where the metabelian eta--invariant sliceness obstruction
is zero, but where the metabelian $L^2$--eta--invariant obstruction is non--zero. 
This shows that more eta--invariants have to
vanish to get zero $L^2$--eta--invariants.

Our main theorem is the following. 

\begin{theorem} \label{thml2eta}
Let $K$ be a knot with zero metabelian eta--invariant ribbonness obstruction,  then
$K$ has zero metabelian $L^2$--eta--invariant sliceness obstruction.
\end{theorem}

The proof of the theorem will be done in the next two sections. 
In \cite{F03} we showed that the converse is not true, i.e. there exists a knot 
with zero metabelian $L^2$--eta--invariant but  non--zero metabelian eta--invariant ribbonness obstruction

\subsection{Approximation of  $L^2$--eta--invariants}
\begin{definition}
We say that $G$ is residually finite it there exists
  a sequence of normal subgroups $G \supset G_1 \supset G_2
\supset \dots$ of finite index $[G:G_i]$ such that $\cap_i G_i=\{1\}$.
We call the sequence $\{G_i\}_{i\geq 1}$ a resolution of $G$.
\end{definition}


If $\varphi: \pi_1(M)\to G$ is a homomorphism to a finite group, 
then define $\eta(M,G)=\eta(M,\a_G)$ where $\a_G:\pi_1(M)\xrightarrow{\varphi} G\to U(\C G)$ 
is the canonical induced unitary
representation given by left multiplication.

\begin{theorem} \label{l2apprthm} \label{finitecoversignature}
Let  $\varphi:\pi_1(M)\to G$ be a homomorphism.
\bn
\item
If $G$ is finite, then
\[ \ba{rcl} \eta(M,G)&=&\sum_{\a \in \hat{R}^{irr}(G)} \dim(\a)\eta(M,_{\a\circ \v}(M) \\
   \eta^{(2)}(M,G)&=&\frac{\eta(M,G)}{|G|} \ea \]
\item
If $G$ is residually finite group then the above equality ``holds in the limit'', i.e.
if $\{G_i\}_{i\geq 1}$ is a resolution of $G$,
then  
\[ \eta^{(2)}(M,G) =\lim_{i\to \infty} \frac{\eta(M,G/G_i)}{|G/G_i|} \]
\en
\end{theorem}

\begin{proof}
The first statement follows immediately from the well-known fact of the representation theory of finite groups that 
 \[ \C G =\sum_{\a \in \hat{R}^{irr}(G)} V_{\a}^{\dim(\a)}\]
The second statement is shown in \cite{A76},
L\"uck and Schick proved the last parts (cf. \cite[remark 1.23]{LS01}).

 \end{proof}

\begin{corollary} \label{corl2equalssum} \label{corl2equalsint}
Let $K$ be a knot, then
\[ 
\ba{rcl} \eta^{(2)}(M_K,\Z/k)&=&\frac{1}{k}\eta(M_K,\Z/k)=\frac{1}{k}\sum_{j=1}^{k}\s_{e^{2\pi ij/k}}(K)\\
  \eta^{(2)}(M_K,\Z)&=&\int_{S^1} \s_{z}(K)\ea  \]
\end{corollary}

This corollary was also proven by Cochran, Orr and Teichner (cf. \cite{COT02}), using a different
approach.

\begin{proof}
The first part is immediate from the decomposition of $\C[\Z/k]$ into one-dimensional
$\C[\Z/k]$--modules.
For the second part consider the sequence $\Z \supset 2! \Z  \supset 3!\Z \supset  4! \Z \supset \dots$,
by theorem \ref{l2apprthm}  and corollary \ref{corl2equalssum}
\[ \eta^{(2)}(M_K,\Z) =\lim_{k\to \infty} \frac{\eta(M_K,\Z/k!)}{k!}=
\lim_{k\to \infty} \frac{ \sum_{j=0}^{k!-1} \s_{e^{2\pi i j/k!}}(K)}{k!}=\int_{S^1} \s_z(K)
 \]
The last equality follows from the fact that $\s_z(K)$ is a step function with only finitely many break points.
\end{proof}

\subsection{Proof of theorem \ref{thml2eta}}

Assume that $K$ has zero metabelian eta--invariant ribbon obstruction.
Let $P$ be a metabolizer such that $\eta(M_K,\a(z,\chi))=0$ for all $\a_{(z,\chi)}\in P_k(\pi_1(M_K))$  with
$\chi(P)\equiv 0$.
Let $P_{\Q}:=P\otimes \Q$, this is a  metabolizer for the rational Blanchfield pairing $\l_{Bl,\Q}$.
 We will show that for any
$x\in P_{\Q}$ 
$\lteta(M_K,\b_x)=0$, where $\b_x$ denotes the  map 
\[ \pi_1(M_K)\to \Z \ltimes H_1(M_K,\L) \to \Z \ltimes H_1(M_K,\Q\L)\xrightarrow{id \times \l_{Bl,\Q}(x,-)}\Z \ltimes
\Q(t)/\Q[t,t^{-1}]\] 
This implies the theorem. 

So let $x\in P_{\Q}$.
Note that $nx\in P$ for some $n\in \N$. The map $\b_{nx}$ factors through 
$\Z \ltimes \Delta_K(t)^{-1}\L/\L$, hence $\b_x$ factors through
$\Z \ltimes n^{-1}\Delta_K(t)^{-1}\L/\L$.
\begin{claim}
 
There exists an isomorphism
\[ \im\{\Z \ltimes n^{-1}\Delta_K(t)^{-1}\L/\L \xrightarrow{}  \Z \ltimes \Q(t)/\Q[t,t^{-1}]\} \to
\Z \ltimes \Delta_K(t)^{-1}\L/\L
\]
\end{claim}
\begin{proof}
 Consider the short exact sequence
\[ 0\to \Delta_K(t)^{-1}\L/\L \to n^{-1}\Delta_K(t)^{-1}\L/\L\to \Delta_K(t)^{-1}\L/n^{-1}\Delta_K(t)^{-1}\L\cong \L/n\to 0\]
since tensoring with $\Q$ is exact and since $\L/n$ is $\Z$--torsion we see that
\[ \im\{ n^{-1}\Delta_K(t)^{-1}\L/\L\to \Delta_K(t)^{-1}\Q \L/\Q\L\}
\cong \mbox{Im}\{\Delta_K(t)^{-1}\L/\L\to \Delta_K(t)^{-1}\Q \L/\Q\L\}\]
But $\Delta_K(t)^{-1}\L/\L\to\Delta_K(t)^{-1}\Q \L/\Q\L\to  \Q(t)/\Q[t,t^{-1}]$ is injective, since
$\Delta_K(t)^{-1}\L/\L$ is
$\Z$--torsion free. This shows that 
\[ \im\{n^{-1}\Delta_K(t)^{-1}\L/\L\to \Q(t)/\Q[t,t^{-1}]\}\cong \Delta_K(t)^{-1}\L/\L \]

 Since all maps preserve the
$\Z$--action the claim follows.
\end{proof}

\begin{lemma} \label{proppp2resfinite}
Let $K$ be a knot, then $\Z \ltimes \Delta_K(t)^{-1}\L/\L$ is residually finite.
\end{lemma}

\begin{proof}
Write
$\Delta_K(t)=a_{2g}t^{2g}+\dots+a_1t+a_0$ with $a_{2g}\ne 0, a_{2g-i}=a_i$. Let
 $p$ be a prime number coprime to  $a_{2g}$. 
Write $H:=\Delta_K(t)^{-1}\L/\L$ and $H_i:=p^i H$.
Then $\{H_i\}_{i \geq 1}$ forms a resolution  for $H$
since 
 there 
exists an embedding $\Delta_K(t)^{-1}\L/\L\cong \L/\Delta_K(t)\L \to \Z[1/a_{2g}]^{2g}$ of
$\Z$-modules.

 Since the $\L$-modules $H/H_i$ are finite there exists for each $i$ a number $k_i$ such 
that $t^{k_i}v= v$  for all $v \in H/H_i$ where $t$ denotes a generator of $\Z$.
 Note that $\Z/k_i\ltimes H/H_i$ and 
the map $ \Z \ltimes H\to \Z \ltimes H/H_i$ are well-defined. 
We can in fact pick $k_i$ with the extra properties that
$ k_i>i$ and $k_i|k_{i+1}$, then it is clear that
the kernels of the maps
\[  \Z \ltimes  \Delta_K(t)^{-1}\L/\L\to \Z/k_i \ltimes H/H_i \]
 define a resolution for $\Z \ltimes \Delta_K(t)^{-1}\L/\L$.

\end{proof}


Let 
\[ G:=\im\{\b_x:\Z \ltimes H_1(M,\L)\to \Z \ltimes n^{-1}\Delta_K(t)^{-1}\L/\L \xrightarrow{}  \Z \ltimes
\Q(t)/\Q[t,t^{-1}]\}\]  Note that $G:=\Z \ltimes H$ for some $H\subset \Delta_K(t)^{-1}\L/\L$.
It follows from the proof of lemma \ref{proppp2resfinite} that we can find $H_i\subset H$ and $k_i$ such that
$H/H_i$ is a $p$-group and such that the kernels $G_i$ of
\[ \Z \ltimes H\to \Z/k_i^{s_i}\ltimes H/H_i\] 
form a resolution
 for any exponents $s_i\in \N$ with
$1\leq s_1\leq s_2\leq \dots$.
We will specify the $s_i$ later. 
Using the fact that in general $\eta^{(2)}(M,\varphi:\pi_1(M)\to J)=\eta^{(2)}(M,\varphi:\pi_1(M)\to \im(J))$
(cf. \cite{COT02})
we get
\[\lteta(M_K,\b_x:\pi_1(M_K)\to \Z \ltimes \Q(t)/\Q[t,t^{-1}])=\lteta(M_K,\b_x:\pi_1(M_K)\to G)\]
The groups $G_i$ are a resolution for $G$, hence by theorem \ref{l2apprthm}
\[ \lteta(M_K,\b_x:\pi_1(M_K)\to G)=\lim_{i\to \infty}\frac{\eta(M_K,G/G_i)}{|G/G_i|}=
\lim_{i\to \infty}\frac{\sum_{\a \in \hat{R}^{irr}(G/G_i)} \dim(\a)\eta(M_K,\a)}{|G/G_i|}
\]
To continue we have to understand the irreducible representations of  $G/G_i\cong \Z/k_i^{s_i} \ltimes H/H_i$.
The proof of the  following lemma is the same as the proof of lemma \ref{lemma1} in \cite{F03}. 

\begin{lemma} \label{lemmarepzk}
Let   $F$ be a finite module over $\L_k:=\Z[t]/(t^k-1)$. Then any irreducible
representation $\Z/k^s \ltimes F\to U(l)$ is  conjugate to
 \[ \a_{(l,z,\chi)}(n,h)=  \a_{(z,\chi)}(n,h):=
    z^n \bp 0& \dots &0&1 \\ 1&\dots &0&0 \\
\vdots &\ddots &&\vdots \\
     0&\dots &1&0 \ep^n \bp \chi(h) &0&\dots &0 \\
 0&\chi(th) &\dots &0 \\
\vdots &&\ddots &\vdots \\ 0&0&\dots &\chi(t^{l-1}h) \ep  \]
for some $z\in S^1$ with $z^k=1$ and $\chi:F\to F/(t^l-1)\to S^1$ a character which does not factor through
$F/(t^r-1)$ for some $r<1$. In particular there are no irreducible representations of dimension greater than $k$.
\end{lemma}

\begin{remark}
Note that $k_i$ is in general a composite number since the order of a $p$--group is always composite. In particular
$\lteta(M_K,\b_x)$ is the limit of eta--invariants which are in general not of prime power dimension. This explains
why the vanishing of the metabelian eta--invariant sliceness obstruction, which involves only
prime power dimensional eta--invariants, does not imply the vanishing of
the
$L^2$--eta--invariant sliceness obstruction.
\end{remark}

This lemma shows that all irreducible representations $\Z\ltimes H_1(M_K,\L)\to G\to G/G_i \cong \Z/k_i^{s_i}\ltimes H/H_i\to
U(l)$ are of the type $\a_{(z,\chi)}$ where $z^{k_i^{s_i}}=1$ and $\chi$ is of prime power order since $H/H_i$ 
is a $p$--group. Furthermore, since $x\in P_{\Q}$ and $P_{\Q}=P_{\Q}^{\perp}$ we have $\chi(P)\equiv 0$.
If the $z$'s had been transcendental our proof would be complete by now since we assumed that
$\eta(M_K,{\a_{(z,\chi)}})=0$ for all $\chi$ of prime power order with $\chi(P)\equiv 0$ and all transcendental $z$.

The next two propositions show that $\eta(M_K,{\a_{(z,\chi)}})=0$ for almost all $z$.
We will see that the non--zero contributions
in $\frac{1}{|G/G_i|}\sum_{\a \in \hat{R}^{irr}(G/G_i)} \dim(\a)\eta(M_K,\a)$ vanish in the limit.


\begin{proposition} \label{corbound}
There exists a number $C$ such that for any $\chi:H_1(M_K,\L)/(t^k-1)\to S^1$ of prime power order the map
\[ \ba{rcl} S^1 &\to &\Z \\
  z &\mapsto &\eta(M_K,{\a(k,z,\chi)}) \ea \]
has at most $Ck$ discontinuities.
 \end{proposition}

 For the proof we need the following lemma.
 
\begin{lemma}{\cite[p. 92]{L94}} \label{propetacont}
Let $M^3$ be a manifold, then for any $r\in \N$ the map
\[ \ba{rcl} \eta_k: R_k(\pi_1(M)) &\to &\R \\
 \a &\mapsto &\eta(M,{\a}) \ea \]
is continuous on  ${\S}_r:=\{ \a \in  R_k(\pi_1(M)) | \sum_{i=0}^3\dim(H_i^{\a}(M,\C^k))=r \}$.
\end{lemma}

Let $J:=\Z \ltimes H_1(M_K,\L)$. Denote the $J$-fold cover of $M_K$ by $\hat{M}$. 
After triangulating $M$ we can view
\[ 0 \to C_3(\hat{M})\xrightarrow{\partial_3}  C_2(\hat{M})
\xrightarrow{\partial_2} C_1(\hat{M})\xrightarrow{\partial_1} C_0(\hat{M})\to 0\]
as a complex of free $\Z J$--modules where $\rk(C_0(\hat{M}))=\rk(C_3(\hat{M}))=1$
and $\rk(C_1(\hat{M}))=\rk(C_2(\hat{M}))=m$ for some $m$. Represent $\partial_2$
 by an $m\times m$-matrix $R$ over $\Z J$.
Then for $\a \in R_k(\pi_1(M_K))$ we get
\[ \det(\a(R))\ne 0 \Rightarrow \a\in {\S}_{2k} \]
since $H_*^{\a}(M,\C^k)=H_*(C_*(\hat{M})\otimes_{\Z J}\C^k)$.

For a character $\chi:H_1(M_K,\L) \to H_1(M_K,\L)/(t^k-1)\to  S^1$ define
\[ S_{k,\chi}:=\{ z\in S^1|\det(\a_{(k,z,\chi)}(R))=0 \} \]

\begin{lemma} \label{propskestimate}
There exists a number $C$ such that $|S_{k,\chi}| \leq Ck$ for all $\chi$
of prime power order.
\end{lemma}

\begin{proof}
Denote by $f:\Z[J]\to \Z[t,t^{-1}]$ the map induced by $(n,v)\mapsto t^n$.
For $g=\sum_{i=n_0}^{n_1} a_it^i, a_{n_0}\ne 0, a_{n_1}\ne 0$ define $\deg(g)=n_1-n_0$.
Let $C:=m\max\{\deg(f(R_{ij}))\}$.
Given a character $\chi$   denote by $z$ a variable, 
then $D(z):=\a_{(z,\chi)}(R)$ is a $km \times km$--matrix over $\C[z,z^{-1}]$. 
It's clear that $\deg(\det(D(z)))\leq \frac{C}{m}km=Ck$, hence either $\det(D(z))\equiv 0$ or there are  at most $Ck$ $z$'s
which are zeroes of $\det(D(z))$.
Letsche \cite[cor. 3.10]{L00} showed that for any $\chi$ of prime power order $S_{k,\chi}$ does not contain
any transcendental number, in particular $\det(D(z))$ is not identically zero.
\end{proof}

This lemma proves proposition \ref{corbound}.

\begin{proposition} \label{propesteta}
For each $k$ there exists $D_k\in \R$ such that
\[ |\eta(M_K,{\a})| \leq D_k \]
for all $\a \in R_l(\pi_1(M_K))$ and all $l\leq k$.
\end{proposition}

\begin{proof}
Let \[ \ti{\S}_r:=\{ \a \in  R_k(\pi_1(M)) | \sum_{i=0}^3\dim(H_i^{\a}(M,\C^k))\geq r \} \]
Levine \cite[p. 92]{L94} shows that these are subvarieties of $R_k(\pi_1(M))$, that $\ti{\S}_N=\varnothing$
for some $N$ and that $\eta_k$ is continuous on $\ti{\S}_r\setminus \ti{\S}_{r+1}$ for all $r$.

We claim that $\eta_k$ is bounded on each $\ti{\S}_r$. 
Note that $\ti{\S}_{r}\setminus \ti{\S}_{r+1}$ has only finitely many components since
$\ti{\S}_{r+1}$ is a subvariety. If $\eta_k$ is not bounded on $\ti{\S}$ then it is therefore 
not bounded on at least one component $C$ of   $\ti{\S}_{r}\setminus \ti{\S}_{r+1}$. 

Since $\pi_1(M_K)$ is finitely generated it follows that $R_l(\pi_1(M_K))$ compact, hence $\bar{C}\subset \ti{\S}_r$ is compact
too. We can therefore find a sequence $p_i \in C$ such that $p_i$ converges to some point $p\in \bar{C}$
and such that $\lim_{i\to \infty} \eta_k(p_i)=\infty$. Since $C$ is path connected and locally path connected
we can find a curve $\g:[0,1]\to C$ such that $\g(1-\frac{1}{2^i})=p_i$. Note that $\g(p[0,1])=[D,\infty)$ for some
$D$. In particular we can find sequences
$q_i$ and $r_i$ in $\ti{\S}_{r}\setminus \ti{\S}_{r+1}$ converging to point $p $ with
$\eta(q_i)=i+\frac{1}{2}$ and $\eta(r_i)=i$. 
But this is a contradiction to the fact, established by Levine \cite[p. 92]{L94}, that 
 $\eta_k \mod \Z:R_k(\pi_1(M))\to \R/\Z$ is continuous. 
 \end{proof}

We are now ready to show that $\eta^{(2)}(M_K,\b_x)=0$ for any $x\in P_{\Q}$ which proves of theorem
\ref{thml2eta}.
Recall that we have to show that
\[  \lim_{i\to \infty}\frac{\eta(M_K,G/G_i)}{|G/G_i|}=0
\] 
We pick $s_i$ with the  extra property $k_i^{s_i-4} \geq D_{k_i}$ for all $i$.
Using lemma \ref{finitecoversignature} we get
\[ |\eta(M_K,G/G_i)|\leq \sum_{\a \in \hat{R}^{irr}(G/G_i)}\dim(\a)|\eta(M,{\a})|\] 
Recall that $G/G_i\cong \Z/k_i^{s_i} \ltimes H/H_i$ and that $H/H_i$ is a $p$-group. By definition of $k_i$
 any character  
actually factors through $(H/H_i)/(t^{k_i}-1)$.
In particular by lemma \ref{lemmarepzk}  there are no irreducible representations of dimension bigger than $k_i$. 
It now follows that the above term is in fact less or equal than
\[ \sum_{j=1}^{k_i}j \sum_{\a \in \hat{R}^{irr}_j(G/G_i)} |\eta(M_K,{\a})|
\leq \sum_{j=1}^{k_i}j
 \sum_{\scriptsize{\chi:(H/H_i)/(t^j-1)\to S^1}} 
 \sum_{z\in S^1,z^{k_i^{s_i}}=1}
|\eta(M_K,{\a(j,z,\chi)})|\]
From  corollary  \ref{corbound} and using that $\eta(M_K,\a_{(z,\chi)})$ for
all transcendental $z$ and all $\chi$ of prime power order with $\chi(P)\equiv 0$, it follows that 
$\eta(M_K,{\a(z,\chi)})=0$ for all but at most $Ck_i$ values of $z$.
Using this observation and using proposition \ref{propesteta} we get that the above term is less or equal 
than
\[ \sum_{j=1}^{k_i} j\sum_{\scriptsize{\chi:(H/H_i)/(t^j-1)\to S^1}} Cj  D_{k_i} \leq k_i^3C|H/H_i|D_{k_i} \]
Therefore
 \[ |\eta^{(2)}(M_K,\b_x)| = \left\arrowvert \lim_{i\to \infty}\frac{\eta(M_K,G/G_i)}{|G/G_i|} \right\arrowvert
\leq \lim_{i\to \infty} \frac{k_i^3CD_{k_i}|H/H_i| }{k_i^{s_i}|H/H_i|}  
= \lim_{i\to \infty} \frac{k_i^3}{k_i^4}\frac{CD_{k_i}}{k_i^{s_i-4}} 
=0
\]
since $\lim_{i\to \infty}k_i=\infty$ and by the choice of $s_i$.
 This concludes the proof of  theorem \ref{thml2eta}.



\end{document}